\documentstyle[12pt]{amsart}
\newtheorem{athm}{Theorem 2.8}

\newtheorem{bthm}{Theorem 2.9}

\newtheorem{cthm}{Theorem 3.5}

\newtheorem{thm}{Theorem}[section]
\newtheorem{lem}[thm]{Lemma}
\begin{document}
\title[Generations for $J_1$ and $J_2$]
{$(p,q,r)$-Generations for the Janko Groups $J_1$ and $J_2$}
\author{Jamshid Moori}
\thanks{supported by a research grant from the University of Natanl}
\maketitle
\begin{center}
Department of Mathematics and Applied Mathematics\\
University of Natal\\
P.O.Box 375, Pietermaritzburg\\
3200 South Africa
\end{center}
\section{Introduction}

Suppose that $G$ is a group which is generated by two elements $a$ and $b$
such that $o(a)=l$ and $o(b)=m$ and their product $ab$ has order $n$. We
say that $G$ is a $(l,m,n)$-generated group and it is obvious that in this
case $G$ must be a quotient of the {\it triangle group} $\Delta(l,m,n)$
given by the presentation
$$ \Delta(l,m,n)=<a,b,c \, | \, a^l = b^m =c^n = abc =1>.$$

If $G$ is a $(l,m,n)$-generated group, then it is also a $(r,s,t)$-generated
group whenever $(r,s,t)$ is a rearrangement of $(l,m,n)$. Thus we may assume
that
\mbox{$l \leq m \leq n$}. It is well-known that $\Delta(l,m,n)$ is finite if
and only
if $ 1/l+1/m+1/n>1$. (see [11] and [4]). Finite $\Delta(l,m,n)$ are
$\Delta(1,n,n)$ the cyclic group of order $n$, $\Delta(2,2,n)$ the dihedral
group of order $2n$, $\Delta(2,3,3)$ the alternating group $A_4$,
$\Delta(2,3,4)$ the symmetric group $S_4$, and $\Delta(2,3,5)$ the
alternating group $A_5$. If $1/l+1/m+1/n = 1$, namely in the cases of
$\Delta(2,3,6)$, $\Delta(2,4,4)$ and $\Delta(3,3,3)$, then the triangle group
is infinite but soluble. In the case when $1/l+1/m+1/n <1$, then the triangle
group is infinite and insoluble. These triangle groups have a remarkable
wealth of interesting finite quotient groups. We encourage reader to consult
[4], [14] and [15] for discussion and background material on the triangle groups.
A $(2,3,7)$-generated group $G$ which give rise to compact Riemann surfaces
of genus greater than 2 with automorphism groups of maximal order, are called
{\it Hurwitz groups} ([7] and [15]).

Woldar in [16] determined that $J_1, J_2, He, Ru, Co_3, HN, Ly$ are Hurwitz,
while $M_{11}, M_{12}, M_{22}, M_{23}, HS, J_3, M_{24}, McL, Suz, O$'$N, Co_1,
Co_2$ are non-Hurwitz. In [8] Kleidman, Parker and Wilson proved that
$F_{23}$ is not Hurwitz, while Linton [9] and Wilson [10] showed that $Th$ and
$F'_{24}$ are Hurwitz. In [17] Woldar completes the problem, except for $B$
and $M$ where the question is unresolved, by proving that $J_4$ and$F_{22}$
are Hurwitz.

If $G$ is a finite group, $C_1, C_2, C_3$ conjugate classes of $G$, and $c$ a
fixed representative of $C_3$, then $\xi_G(C_1,C_2,C_3)$ denotes the {\it
structure constant} of the group algebra $ Z({\Bbb C}[G])$ which is equal to
the number of distinct ordered pairs $(a,b)$ with $a \in C_1$, $b \in C_2$ and
$ab=c$. Using the character table of $G$, the number $\xi_G(C_1,C_2,C_3)$
can be calculated by the formula
$$\xi_G(C_1,C_2,C_3)=\frac{|C_1|.|C_2|}{|G|}\sum_{i=1}^{k} \frac{\chi_i(a)
\chi_i(b) \overline{\chi_i(c)}}{\chi_i(1)}$$
where $\chi_1, \chi_2, \ldots ,\chi_k$ are the irreducible complex characters
of $G$. Let $\xi_G^*(C_1,C_2,C_3)$ denote the number of distinct ordered pairs
$(a,b)$ with $a \in C_1$, $b \in C_2$, $ab=c$ and $G=<a,b>$. Obviously $G$ is a
$(l,m,n)$-generated group if and only if there exist conjugacy classes $C_1,
C_2, C_3$ with representatives $a, b, c$, respectively such that $o(a)=l,
o(b)=m$ and $o(c)=n$, for which $\xi_G^*(C_1,C_2,C_3)>0.$ In this case we
say that $G$ is $(C_1,C_2,C_3)$-generated. If $H$ is a subgroup of $G$
containing $c$ and $B$ is a conjugacy class of $H$ such that $c \in B$, then
$\sigma_H(C_1,C_2,B)$ denotes the number of distinct ordered pairs $(a,b)$
such that $a \in C_1$, $b \in C_2$, $ab=c$ and $<a,b> \leq H$.

It is well-known that every finite non-abelian simple group can be generated
by two suitable elements. (For details see [5]). If $G$ is non-abelian finite
simple group and $l,m,n$ are divisors of $|G|$ such that $1/l+1/m+1/n<1$,
then the following question arises: Is $G$ a $(l,m,n)$-generated group? We
are a great distance away from a complete answer to the question, however we
are aiming at a partial answer to it by considering the case when $l=p,m=q$
and $n=r$ are primes and the group $G$ is isomorphic to one of the sporadic
simple groups.

In the present paper we investigate $(p,q,r)$-generations for the Jank groups
$J_1$ and $J_2$ where $p,q$ and $r$ are distinct primes satisfying \mbox {$p
<q<r$}. We prove the following results.
\begin{athm}
The group $J_1$ is $(p,q,r)$-generated for $p,q,r \in \{2,3,5,7,11,19\}$
with $p<q<r$, except when $(p,q,r)=(2,3,5)$.
\end{athm}
\vspace{-.7 cm}
\begin{bthm}
The group $J_1$ is generated by three involutions $a,b,c \in 2A$ such
that $abc \in 11A$.
\end{bthm}
\vspace{-.7 cm}
\begin{cthm}

\begin{enumerate}
\item[(a)] The group $J_2$ is $(2B,3B,7A)$-, $(2A,5C,7A)$-, $(2B,5A,7A)$-,
$(2B,5C,7A)$-,\\
$(3A,5C,7A)$-, $(3B,5A,7A)$-, and $(3B,5C,7A)$-generated.
\item[(b)] The group $J_2$ is not $(2A,3A,7A)$-, $(2A,3B,7A)$-, $
(2B,3A,7A)$-,
\\
$(2A,5A,7A)$-, or $(3A,5A,7A)$-generated.
\end{enumerate}
\end{cthm}

The remaining cases , namely when $p=q=r$, $p \neq q=r$ or $p=q \neq r$
will be treated separately in another paper by the author.

For the description of the conjugacy classes, the character
tables and information
on the maximal subgroups readers are referred to ATLAS [3]. Computations were
performed with the aid of CAYLEY and GAP (see [1] and [13]) running on a
\mbox{ SUN GX2} computer.
\subsection{Acknowledgments}
I am very grateful to the referee for his valuable and constructive
remarks, and to Jonathan Hall for several of the arguments in this paper.
I also thank Andrew Woldar for his private communication on the \mbox {Lemma
2.1}.

\section{$(p,q,r)$-generations for $J_1$}

We list bellow the structure constants for the Janko group $J_1$.\\

\begin{center}
\begin{tabular}{c|rrr}
$\xi_{J_1}(C_1,C_2,C_3)$ & $C_1$ & $C_2$ & $C_3$\\
\hline
$49$ & $2A$ & $3A$ & $7A$\\
$55$ & $2A$ & $3A$ & $11A$\\
$38$ & $2A$ & $3A$ & $19A$\\
$49$ & $2A$ & $5A$ & $7A$\\
$44$ & $2A$ & $5A$ & $11A$\\
$57$ & $2A$ & $5A$ & $19A$\\
$209$ & $2A$ & $7A$ & $11A$\\
$209$ & $2A$ & $7A$ & $19A$\\
$133$ & $2A$ & $11A$ & $19A$\\
$189$ & $3A$ & $5A$ & $7A$\\
$198$ & $3A$ & $5A$ & $11A$\\
$171$ & $3A$ & $5A$ & $19A$\\
$858$ & $3A$ & $7A$ & $11A$\\
$836$ & $3A$ & $7A$ & $19A$\\
$494$ & $3A$ & $11A$ & $19A$\\
$858$ & $5A$ & $7A$ & $11A$\\
$836$ & $5A$ & $7A$ & $19A$\\
$513$ & $5A$ & $11A$ & $19A$\\
$2299$ & $7A$ & $11A$ & $19A$\\

\end{tabular}
\end{center}

We would also like to mention here that if $G$ is a $(p,q,r)$-generated group
where $p,q$ and $r$ are distinct primes, then $G$ has no soluble quotient.
This is an obvious consequence of the following elementary result:

\begin{lem}
Let $G$ be a $(l,m,n)$-generated group with $l,m,n$ pairwise co-prime.
Then $G$ has no soluble quotient.
\end{lem}

\begin{pf}
Assume that $G$ is generated by $a,b$ such that $o(a)=l, o(b)=m$ and \mbox {
$o(ab)=n$.} If $G$ has soluble quotient, then it also has a non-trivial abelian quotient
$A$. Let $\bar{a}, \bar{b}$ and $\overline{ab}$ be the respective images of
$a,b$ and $ab$ in $A$. Then $o(ab)$ divides $lm$ and $n$. Since $lm$ and $n$
are relatively prime, $\overline{ab}= \bar{e}$, where $e$ is the identity of
$G$. This implies that $\bar{a}= (\bar{b})^{-1}$. So if $\bar{b} \neq \bar{e}$
then $o(b)$ divides both $l$ and $m$, a contradiction to co-primeness. But if
$\bar{b} = \bar{e}$, then $\bar{a} = \bar{e}$. Hence $A=<\bar{a},\bar{b}>=\{
\bar{e}\}$, contradicting the non-triviality of $A$. This completes the proof
of the lemma.
\end{pf}

We now deal with the $(p,q,r)$-generations of $J_1$.

\begin{lem}
The group $J_1$ is $(2,3,7)$-generated.
\end{lem}
\begin{pf}
The fact that $J_1$ is $(2A,3A,7A)$-generated is given by Woldar in [16].
\end{pf}

\begin{lem}
The group $J_1$ is $(2,3,11)$-generated.
\end{lem}
\begin{pf}
The only non-soluble maximal subgroups of $J_1$ with order divisible by 11
are isomorphic to the group $L_2(11)$. If $H \leq J_1$ with $H \cong L_2(11)$,
then
$$\sigma_H(2A,3A,11A) = \xi_H(2A,3A,11A) = 11. $$
We also observe that $N_H(11A) = 11 \! : \! 5$ and $N_{J_1}(11A) = (11 \! :
\! 5) \! : \! 2$. Let
$2=<\alpha>$. Then a fixed element in $11A$ is contained in $H$ and $H^g$ for
$g \in {J_1-H}$ if and only if $g=h^{-1}\beta\alpha$ for some $h \in H$ and
$\beta \in N_H(11A)$. Thus a fixed element of order $11$ in $J_1$ is contained
in precisely two copies of $H$, namely $H$ and $H^\alpha$. Hence $H$ and its
conjugates contribute $11 \times 2 =22$ to the $\xi_{J_1}(2A,3A,11A)$. Since
$\xi_{J_1}(2A,3A,11A)=55>22$, $\xi^*_{J_1}(2A,3A,11A) \geq 33$. Therefore $J_1$
is $(2A,3A,11A)$-generated.
\end{pf}

\begin{lem}
The group $J_1$ is $(2,3,19)$-generated.
\end{lem}
\begin{pf}
The only maximal subgroup of $J_1$ with order divisible by $19$, up to
isomorphism, is the group $19 \! : \! 6$. Since a $(2,3,19)$-generated group cannot
have a soluble quotient, $J_1$ posses no proper $(2,3,19)$-generated
subgroup. Now the result follows from the fact that $\xi_{J_1}(2A,3A,19A)=
\xi^*_{J_1}(2A,3A,19A) = 38$.
\end{pf}

\begin{lem}
The group $J_1$ is $(2,5,11)$-generated.
\end{lem}
\begin{pf}
Up to isomorphism, $L_2(11)$ is the only non-soluble maximal subgroup of
$J_1$ with order divisible by $11 \times 10=110$. We also have
$\sigma_{L_2(11)}(2A,5A,11A)=\xi_{L_2(11)}(2A,5A,11A)=11$, and as in
Lemma 2.3, $L_2(11)$ and it's conjugates contribute $11 \times 2 =22$ to
the $\xi_{J_1}(2A,5A,11A)$. Since $\xi_{J_1}(2A,5A,11A)=55$, result follows.
\end{pf}

\begin{lem}
The group $J_1$ is $(2,5,7)$-, $(2,5,19)$-, $(2,7,11)$-, $(2,7,19)$-,\\
$(2,11,19)$-, $(5,7,11)$-, $(5,7,19)$-, $(5,11,19)$-, $(3,5,7)$-,
$(3,5,19)$-,\\
$(3,7,11)$-, $(3,7,19)$-, $(3,11,19)$-, and $(7,11,19)$-generated.
\end{lem}
\begin{pf}
The group $J_1$ has no proper subgroups of order divisible by $2 \times 5
\times 7$, $2 \times 5 \times 19$, $2 \times 7 \times 11$, $2 \times 7 \times
19$, $2 \times 11 \times 19$, $5 \times 7 \times 11$, $5 \times 7 \times 19$
, $5 \times 11 \times 19$, $3 \times 5 \times 7$, $3 \times 5 \times 19$,
$3 \times 7 \times 11$, $3 \times 7 \times 19$, $3 \times 11 \times 19$, and
 $7 \times 11 \times 19$ respectively. Since $\xi_{J_1}(C_1, C_2, C_3)$ for
the corresponding conjugacy classes are positive, the result follows.
\end{pf}

\begin{lem}
The group $J_1$ is $(3,5,11)$-generated.
\end{lem}
\begin{pf}
Up to isomorphism, $L_2(11)$ is the only maximal subgroup of $J_1$ with order
divisible by $3 \times 5 \times 11$. Now $\sigma_{L_2(11)}(3A,5A,11A)=
\xi_{L_2(11)}(3A,5A,11A)=22$ implies that $L_2(11)$ and its conjugates
contribute $22 \times 2 =44$ to the $\xi_{J_1}(3A,5A,11A)$. Since
$\xi_{J_1}(3A,5A,11A)=198$, the result follows.
\end{pf}

\begin{thm}
The group $J_1$ is $(p,q,r)$-generated for $p,q,r \in \{2,3,5,7,11,19\}$ with
$p<q<r$, except when $(p,q,r)=(2,3,5)$.
\end{thm}
\begin{pf}
This follows from the Lemmas 2.2, 2.3, 2.4, 2.5, 2.6, 2.7, 2.8
and the fact $\Delta(2,3,5) \cong A_5$.
\end{pf}

If $C_1, C_2, C_3, C_4$ are conjugacy classes of a group $G$, and $d$ a fixed
representative of $C_4$ then $\xi_G(C_1,C_2,C_3,C_4)$ denotes the number of
distinct triples $(a,b,c)$ with $a \in C_1$, $b \in C_2$ and $c \in C_3$
such that $abc=d$. This number is computed using the formula
$$\xi_G(C_1,C_2,C_3,C_4)= \frac{|C_1||C_2||C_3|}{|G|}\sum_{i=1}^{k}\frac{\chi_i(a) \chi_i(b) \chi_i(c) \overline{\chi_i(d)}}{(\chi_i(1))^2} $$
when $\chi_1, \chi_2, \ldots , \chi_k$ are the irreducible complex characters of
$G$. We use this fact to prove the following theorem.

\begin{thm}
The group $J_1$ is generated by three involutions $a,b,c \in 2A$ such that
$abc \in 11A$.
\end{thm}
\begin{pf}
Using the character table of $J_1$ we have $\xi_{J_1}(2A,2A,2A,11A)=17908$.
In $J_1$ we have two maximal subgroups, up to isomorphism, with order
divisible by $11$, namely $11 \! : \! 10$ and $L_2(11)$. Here $N_{J_1}(11A)=
11 \! : \! 10$.
We also have
$$\sigma_{L_2(11)}(2A,2A,2A,11A)=\xi_{L_2(11)}(2A,2A,2A,11A)=242. $$
A fixed element of order 11 in $J_1$ lies in two conjugates of $L_2(11)$.
Hence $L_2(11)$ contributes $242 \times 2=484$ to the number
$\xi_{J_1}(2A,2A,2A,11A)$. It is routine to compute the character table of
the group $11 \! : \! 10$. We represent a part of this character table, giving the
values of irreducible characters on the classes $1A$, $2A$, and $11A$ in
Table 1.
\begin{table}[htb]
\begin{center}
\begin{tabular}{l|rrr}
Centralizer & 110 & 10 & 11\\
\hline
Class & $1A$ & $2A$ & $11A$\\
\hline
$\chi_1$ & 1 & 1 & 1\\
$\chi_2$ & 1 & -1 & 1\\
$\chi_3$ & 1 & 1 & 1\\
$\chi_4$ & 1 & -1 & 1\\
$\chi_5$ & 1 & 1 & 1\\
$\chi_6$ & 1 & -1 & 1\\
$\chi_7$ & 1 & 1 & 1\\
$\chi_8$ & 1 & -1 & 1\\
$\chi_9$ & 1 & 1 & 1\\
$\chi_{10}$ & 1 & -1 & 1\\
$\chi_{11}$ & 10 & 0 & 1\\
\end{tabular}
\end{center}
\caption{Partial character table of $11 \! : \! 10$}
\end{table}

Using Table 1 we have
\begin{align*}
\sigma_{11:10}(2A,2A,2A,11A)
&= \frac{|2A|^3}{11 \times 10} \sum_{i=1}^{k} \frac{(\chi_i(2A))^3 \chi_i(11A)}
{(\chi_i(1))^2}\\
&= \frac{|2A|^3}{11 \times 10} \sum_{i=1}^{10} (-1)^{i+1}\\
&= 0.
\end{align*}
Hence $11\!:\!10$ does not contribute to the number $\xi_{J_1}(2A,2A,2A,11A)$.
Since $\xi_{J_1}(2A,2A,2A,11A) = 17908>484$, the group $J_1$ is $(2A,2A,2A,
11A)$-generated. This completes the theorem.
\end{pf}

\section{$(p,q,r)$-generations for $J_2$}

We list below the structure constants for the group $J_2$.\\

\begin{center}
\begin{tabular}{c|rrr}
$\xi_{J_2}(C_1,C_2,C_3)$ & $C_1$ & $C_2$ & $C_3$\\
\hline
0 & $2A$ & $3A$ & $7A$\\
7 & $2A$ & $3B$ & $7A$\\
0 & $2B$ & $3A$ & $7A$\\
70 & $2B$ & $3B$ & $7A$\\
0 & $2A$ & $5A$ & $7A$\\
7 & $2A$ & $5C$ & $7A$\\
7 & $2B$ & $5A$ & $7A$\\
49 & $2B$ & $5C$ & $7A$\\
0 & $3A$ & $5A$ & $7A$\\
14 & $3A$ & $5C$ & $7A$\\
56 & $3B$ & $5A$ & $7A$\\
343 & $3B$ & $5C$ & $7A$\\

\end{tabular}
\end{center}

\begin{lem}
The group $J_2$ is not $(2A,3A,7A)$-, $(2B,3A,7A)$-, $(2A,5A,7A)$-,
or $(3A,5A,7A)$-generated.
\end{lem}
\begin{pf}
This follows trivially since $\xi_{J_2}(C_1,C_2,C_3)= 0$ for the
corresponding conjugacy classes of $J_2$.
\end{pf}

\begin{lem}
The group $J_2$ is $(2A,5C,7A)$-, $(2B,5A,7A)$-, $(2B,5C,7A)$-,
\\
$(3A,5C,7A)$-, $(3B,5A,7A)$-, $(3B,5C,7A)$-generated.
\end{lem}
\begin{pf}
The only maximal subgroups of $J_2$ with order divisible by 7, up to
isomorphism, are the groups $U_3(3)$ and $L_3(2) \! : \! 2$. Since $|U_3(3)|
 = 2^5.3^2.7$ and $|L_3(2) \! : \! 2|= 2^4.~3.~7$ and $5$ does not divide
 neither $|U_3(3)|$ nor $|L_2(3) \! : \! 2|$, we can say that no proper
subgroup of $J_2$ is $(p,5,7)$-generated when $p \in \{2,3\}$. As
$$\xi_{J_2}(2A,5C,7A)= \xi_{J_2}(2B,5A,7A)=7 \quad ,
\quad \xi_{J_2}(2B,5C,7A)=49$$
$$\xi_{J_2}(3A,5C,7A)=14 \quad , \quad \xi_{J_2}(3B,5A,7A)=56 \quad , \quad
\xi_{J_2}(3B,5C,7A)=343$$
the result will immediately follow.
\end{pf}

\begin{lem}
The group $J_2$ is $(2B,3B,7A)$-generated.
\end{lem}
\begin{pf}
This is given in [6].
\end{pf}

\begin{lem}
The group $J_2$ is not $(2A,3B,7A)$-generated.
\end{lem}
\begin{pf}
Here we use a theorem of Ree (see [12] and [2]) which states: {\it If $G$
is a transitive permutation group generated by permutations $g_1,g_2,\ldots
 ,g_s$ acting on a set of $n$ elements such that $g_1g_2 \cdots g_s$
is the identity permutation, and if generator $g_i$ has exactly $c_i$ cycles
for $1 \leq i \leq s$, then \mbox{$c_1+c_2+ \cdots +c_s \leq (s-2)n + 2$.}}

The group $J_2$ acts as a transitive rank-$3$ group on a set $X$ of $100$
elements. The point stabilizer is isomorphic to the group $U_3(3)$ with
orbits of lengths $1$, $36$, and $63$. If we denote the character of the
permutation representation of $J_2$ on the set $X$ by $\chi$, then by
referring to the character table of $J_2$ we obtain \mbox{$\chi =
1a+36a+63a$} where $ma$ is the first irreducible character of degree $m$ in
the ATLAS character table of $J_2$. It is easy to verify that $\chi(2A)=20$,
$\chi(3B)=4$, and $\chi(7A)=2$. This implies that, in the action of $J_2$
on the set $X$, the elements in $2A$, $3B$ and $7A$ induce permutations with
cycle types $1^{20}2^{40}$, $1^43^{32}$ and $1^27^{14}$ respectively. Now
$$c_1+c_2+c_3 = (20+40)+(4+32)+(2+14)= 112$$ and $$(s-2)n +2=(3-2) \times 100
+2 = 102$$ imply that $c_1+c_2+c_3>(s-2)n+2$. This contradicts Ree's theorem
stated above. Hence $J_2$ is not $(2A,3B,7A)$-generated.
\end{pf}

\begin{thm}
For the Janko group $J_2$ we have the following.
\begin{enumerate}
\item[(a)] $J_2$ is $(2B,3B,7A)$-, $(2A,5C,7A)$-, $(2B,5A,7A)$-,
$(2B,5C,7A)$-,\\ $(3A,5C,7A)$-, $(3B,5A,7A)$-,
and $(3B,5C,7A)$-generated.
\item[(b)] $J_2$ is not $(2A,3A,7A)$-, $(2A,3B,7A)$-,
$(2B,3A,7A)$-, $(2A,5A,7A)$-,\\ or $(3A,5A,7A)$-generated.
\end{enumerate}
\end{thm}
\begin{pf}
This follows from the Lemmas 3.1, 3.2, 3.3 and 3.4.
\end{pf}

\end{document}